\documentclass{article}
\usepackage{latexsym}
\usepackage{amsmath}
\usepackage{amssymb}

\usepackage{theorem}
\usepackage{epsfig}

\newtheorem{theorem}{Theorem}[section] 
\newtheorem{lemma}[theorem]{Lemma}

\newtheorem{corollary}[theorem]{Corollary}
\newtheorem{conjecture}[theorem]{Conjecture}

{\theorembodyfont{\rmfamily}
\theoremstyle{plain}
\newtheorem{definition}[theorem]{Definition}

\newtheorem{remark}[theorem]{Remark}
\newtheorem{ack}{Acknowledgements}
}

\newcommand{\A}{\ensuremath{\mathcal A}}
\renewcommand{\S}{\ensuremath{\mathcal S}}
\newcommand{\C}{\ensuremath{\mathbb{C}}}

\newcommand{\Z}{\ensuremath{\mathbb{Z}}}
\newcommand{\La}{\ensuremath{\Lambda}}

\newcommand{\e}{\ensuremath{\epsilon}} 
\renewcommand{\bar}[1]{\ensuremath{\overline{#1}}}
\newcommand{\rk}{\ensuremath{\mbox{\it rk}}}
\newcommand{\qed}{\hfill \mbox{$\Box$}\medskip}
\newenvironment{proof}{\noindent {\it proof:}}{\qed \par}
\newenvironment{proofof}[1]{\noindent {\it proof of #1:}}{\qed \par}

\title{Orlik-Solomon algebras and Tutte polynomials\thanks{research conducted under an NSF {\em
Research Experiences for Undergraduates} grant.}} 

\author{Carrie J. Eschenbrenner 
  \and Michael J. Falk }

\date{}

\begin{document}
 
\maketitle

\begin{abstract}

The $OS$ algebra $A$ of a matroid $M$ is a graded algebra related to the Whitney homology of the lattice of flats of $M$. In case $M$ is the underlying matroid of a hyperplane arrangement \A\ in $\C^r$, $A$ is isomorphic to the cohomology algebra of the complement $\C^r\setminus \bigcup \A.$ Few examples are known of pairs of arrangements with non-isomorphic matroids but isomorphic $OS$ algebras. In all known examples, the Tutte polynomials are identical, and the complements are homotopy equivalent but not homeomorphic.

We construct, for any given simple matroid $M_0$, a pair of infinite families of matroids $M_n$ and $M'_n$, $n\geq 1$, each containing $M_0$ as a submatroid, in which corresponding pairs have isomorphic $OS$ algebras. If the seed matroid $M_0$ is connected, then $M_n$ and $M'_n$ have different Tutte polynomials. As a consequence of the construction,  we obtain, for any $m$, $m$ different matroids with isomorphic $OS$ algebras. Suppose one is given a pair of central complex hyperplane arrangements $\A_0$ and $\A_1$. Let $\S$ denote the arrangement consisting of the hyperplane $\{0\}$ in $\C^1$. We define the parallel connection $P(\A_0,\A_1)$, an arrangement realizing the parallel connection of the underlying matroids, and show that the direct sums $\A_0 \oplus \A_1$ and $\S\oplus P(\A_0,\A_1)$ have diffeomorphic complements. 

\end{abstract}

\begin{section}{Introduction}

Let $M$ be a simple matroid with ground set $E$. Associated with $M$ is a
graded-commutative
algebra $A(M)$ called the {\em Orlik-Solomon} ($OS$) algebra of $M$. Briefly,
$A(M)$ is the quotient of the free exterior algebra $\Lambda(E)$ on $E$ by
the ideal generated by ``boundaries" of circuits in $M$. If
$\A$
is an arrangement in $\C^r$ realizing the matroid $M$, then $A(M)$ is
isomorphic
to the cohomology algebra of the complement $C(\A)=\C^r \setminus
\bigcup
\A$. So in the attempt to classify  homotopy types of complex
hyperplane complements
one is led to study graded algebra isomorphisms of $OS$ algebras.

The structure of $A(M)$ as a graded vector space is determined uniquely
by the characteristic polynomial $\chi_M(t)$ of $M$. In most cases,
even for matroids having the same characteristic polynomial, the $OS$
algebras can be distinguished using more delicate invariants of the
multiplicative structure \cite{F6,F7}.
In \cite{F3}, however, two infinite families of rank three matroids are
constructed in which corresponding pairs have isomorphic $OS$ algebras, generalizing a result of L.~Rose and H.~Terao \cite[Example 3.77]{OT}. 

The Tutte polynomial $T_M(x,y)$ is an invariant of $M$ that specializes to
$\chi(t)$ under the substitution $x=1-t, y=0$. In the
examples referred to above, the associated matroids have identical Tutte
polynomials. Furthermore, in \cite{F7} it is shown that, under a fairly
weak hypothesis which is satisfied in all known cases, the Tutte polynomial
of a rank three matroid $M$ can be reconstructed from $A(M)$. It is natural
to conjecture that $A(M)$ determines $T_M(x,y)$ in general. The purpose of
this paper is to show that, without additional hypotheses, counterexamples to 
this conjecture abound.
Here is our main result.

\begin{theorem} Let $M_0$ be an arbitrary connected matroid without loops or multiple points. 
Then for each positive
integer $n\geq 3$, there exist matroids $M_n$ and $M'_n$ of rank $\rk(M_0)+n-1$
satisfying
\begin{enumerate}
\item $M_0$ is a submatroid of $M_n$ and $M'_n$.
\item $A(M_n)$ is isomorphic to $A(M'_n)$ as a graded algebra.
\item $T_{M_n}(x,y)\not = T_{M'_n}(x,y)$.
\end{enumerate}
\label{main}
\end{theorem}
In the other direction, we find several examples in \cite{F7} of matroids with
the same Tutte polynomials and non-isomorphic $OS$ algebras.

The matroid $M_n$ of the theorem is simply the direct sum of $M_0$ with
the polygon matroid $C_n$ of the $n$-cycle. The matroid $M'_n$ can be taken to
be the direct sum of an isthmus with {\em any} parallel connection of $M_0$
and $C_n$. 
Thus, by careful choice of $M_0$,
we obtain the following corollary.

\begin{corollary}Given any positive integer $m$, there exist $m$
nonisomorphic simple matroids with isomorphic $OS$ algebras.
\label{cor}
\end{corollary}
Note that the matroids $M_n$ and $M'_n$ have rank greater than three,
and neither is connected. So it remains possible that $A(M)$ determines
$T_M(x,y)$ for matroids of rank three, or for connected matroids.

The arrangements constructed in \cite{F3} were shown to have homotopy equivalent
complements, and the isomorphism of $OS$
algebras is a corollary.
In the last section we prove a far more general result in the high rank setting
of the present work. We define the parallel connection $P(\A_0,\A_1)$ of 
two arrangements in Section 4,
as a natural realization of the parallel connection of the underlying matroids.
The direct sum $\A_0\oplus \A_1$, denoted by $\A_0\coprod \A_1$ in
\cite{OT}, realizes the direct sum of
the underlying matroids.

\begin{theorem} Let $\A_0$ and $\A_1$ denote arbitrary arrangements. Let
\S\ denote the unique nonempty central arrangement of rank 1. Then
$\A_0 \oplus \A_1$
and $\S \oplus P(\A_0,\A_1)$ have diffeomorphic complements.
\label{homo}
\end{theorem}

The examples of \cite{F3} are generic sections
of the arrangements described in Theorem \ref{homo}, 
with $\A_0$ and $\A_1$  
of rank two. The fact that their fundamental groups are isomorphic then follows
immediately from \ref{homo} by the Lefshetz hyperplane theorem. 
For these particular arrangements,
the complements are homotopy equivalent. It is possible that for more
more general $\A_0$ and $\A_1$, this construction could yield rank-three arrangements
whose complements have isomorphic fundamental groups but are not homotopy
equivalent. This phenomenon has not been seen before, and would be of considerable interest.
Also worthy of note is the result of \cite{JY1} that, for arrangements of rank
three, the diffeomorphism type of the complement determines the underlying matroid.
Theorem \ref{homo} demonstrates that this result is false in ranks greater than
three.

The formulation of Theorem \ref{homo} affords a really easy
proof, providing an alternative for the proof of Theorem
\ref{main}(ii), in case $M_0$ is a realizable matroid. The proof is based on a
simple and well-known relation \cite{B,OT} between the topology of the complement of a
central arrangement in $\C^r$ and that of its projective image, 
which coincides
with the complement of an affine arrangement in $\C^{r-1}$, called the  ``decone" of \A.
The proof of Theorem \ref{homo} demonstrates that all the known cases where topological invariants
coincide even while underlying matroids differ are consequences of this 
fundamental principle.

Here is an outline of the proof of Theorem \ref{main}. 
Once the matroids $M_n$ and $M'_n$ are
constructed in the next section, we define a map at the exterior algebra level
which is easily seen to be an  isomorphism. We show that this map carries
relations to relations, hence induces a well-defined map $\phi$ 
of $OS$ algebras.
This map is automatically surjective. Then in Section 4 we compute the 
Tutte polynomials
of $M_n$ and $M'_n$. These are shown to be unequal provided $M_0$ is
connected, but they coincide upon specialization to $y=0$. Thus $M_n$ and
$M'_n$ have identical characteristic polynomials.
It follows that the $OS$ algebras have the same dimension in each degree,
so that $\phi$ must be injective.
In the final section we prove Corollary \ref{cor} and Theorem \ref{homo},
and close with a few comments and a conjecture.

\end{section}

\begin{section}{The construction}

We refer the reader to \cite{Wh1,Ox} for background material on 
matroid theory and Tutte polynomials,
and to \cite{OT} for more information on arrangements and $OS$
algebras.

Let $C_n$ be the polygon matroid of the $n$-cycle. Thus $C_n$ is a matroid of
rank $n-1$ on $n$ points, with one circuit, of size $n$.  This matroid 
is realized by any arrangement
$\A_n$
of $n$ hyperplanes in general position in $\C^{n-1}$. The ground set of
$C_n$ will be taken to be $[n]:=\{1,\ldots ,n\}$ throughout the paper.

Fix a simple matroid $M_0$ with ground set $E_0$ disjoint from $[n]$. Thus $M_0$ has no loops or multiple points.
Let $M_n=C_n \oplus M_0$. 
So the circuits of
$M_n$ are those of $M_0$ together with the unique circuit of $C_n$. If $\A_0$
is an arrangement realizing $M_0$ in $\C^r$, then $M_n$ is realized 
by the direct sum
of $\A_n$ and $\A_0$ in $\C^{r+n-1}$, denoted $\A_n \coprod \A_0$ in \cite{OT}.

Now fix $\e_0\in E_0$. Let $P^n_{\e_0}$ denote the parallel connection
$P(C_n,M_0)$ of $C_n$
with $M_0$ along $\e_0$. Loosely speaking, $P^n_{\e_0}$ is the freest matroid obtained from $C_n$
and $M_0$ by identifying $\e_0$ with the point $1$ of $C_n$. 
Here is a precise definition. Define an equivalence relation on $E:=[n]\cup
E_0$ so that $\{1,\e_0\}$ is the only nontrivial equivalence class. Denote
the class of any $p\in E$ by $\bar{p}$. For $X\subset E$ let $\bar{X}$ be
the set of classes of elements of $X$. Then $P^n_{\e_0}$ is the matroid 
on the set $\bar{E}$ whose set of circuits is 
\begin{multline*}
{\mathcal C}=\{\bar{C} \ |
\ C \ \text{is a circuit of} \ C_n \ \text{or} \ M_0\}  \\ \cup
\{\bar{C-1}\ \cup \ \bar{C'-\e_0} \ | \ 1 \in C \ \text{a circuit of} \ C_n \\
\text{and} \ \e_0 \in C' \ \text{a circuit of} \ M_0\}.
\end{multline*}

Let $S$ denote an isthmus, that is, the matroid of rank one on a
single point, which point will be denoted $p$. Finally, let $M'_n$ be the
direct sum $P^n_{\e_0}\oplus S$.

These two matroids $M_n$ and $M'_n$ are most easily understood in terms of
graphs. If $M_0$ is a graphic matroid, then $M_n$ is the polygon matroid of the
union  (with or without a vertex in common) 
of the corresponding graph $G$ with the $n$-cycle. The
parallel connection $P^n_{\e_0}$ is the matroid of the graph obtained by
attaching a path of length $n-1$ to the vertices of an edge $\e_0$ of $G$,
and $M'_n$ is then obtained by throwing in a pendant
edge $p$. These graphs are illustrated in Figure \ref{graphs}, with $n=6$.

%\psdraft
\begin{figure}
\begin{center}
\epsfig{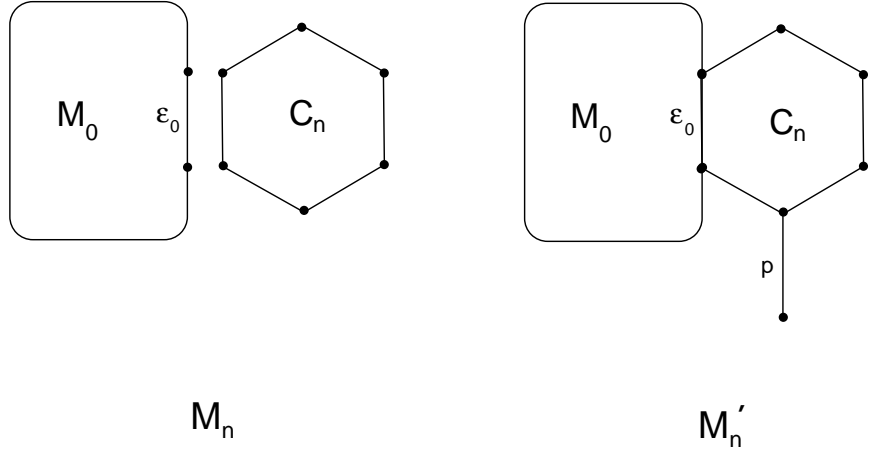}
\caption{The construction.} \label{graphs}
\end{center}
\end{figure}
%\psfull

\end{section}

\begin{section}{An algebra homomorphism}

We proceed to define the $OS$ algebra of a matroid.
Let $M$ be a simple matroid with ground set $E$. Let $\La =\La(E)$ be the
free exterior
algebra 
generated by degree one elements $e_i$ for $i\in E$. 
The results of this paper will hold
for coefficients in
any commutative ring. Define $\partial: \La(E) \to
\La(E)$ by $$\partial(e_1\cdots e_k)=\sum_{i=1}^k (-1)^{i-1}e_1\cdots
\hat{e_i}\cdots e_k,$$ and extending to a linear map.
Let $I=I(M)$ be the ideal of $\La(E)$ generated by $$\{\partial(e_1\cdots
e_k) \ | \ \{e_1,\ldots e_k\} \ \text{is a circuit of} \  M\}.$$
\begin{definition} The {\em $OS$ algebra} $A(M)$ of $M$ is the quotient
$\La(E)/I(M)$.
\end{definition}
Since \La\ is graded and  $I$ is generated by homogeneous elements, 
$A(M)$ is a graded algebra.

The definition of $A(M)$ is motivated by differential topology. Suppose
$\A=\{H_1,\ldots, H_n\}$ is an arrangement of hyperplanes in $\C^r$
realizing the matroid $M$. Let $C(\A)=\C^r - \bigcup_{i=1}^n H_i$.
Extending work of V.I.~Arnol'd and E.~Brieskorn,  P.~Orlik and L.~Solomon
proved the following theorem \cite{OS1}.

\begin{theorem} The cohomology algebra $H^*(C(\A),\C)$ of the complement
$C(\A)$ is isomorphic to $A(M)$.
\label{OS}
\end{theorem}

We now specialize to the examples constructed in the last section. For
simplicity we suppress much of the notation. Consider the integer $n\geq 3$,
the matroid $M_0$, and the point $\e_0$ to be fixed once and for all. 
Unprimed symbols $M,\La,I,A$ will refer to the matroid $M_n$, and primed
symbols $M',\La',I',A'$ refer to $M'_n$.

Recall the ground sets of $M$ and
$M'$ are $E=[n]\, \cup \, E_0$ and $\bar{E}\, \cup\, \{p\}$ respectively. The
generator of $A'$ corresponding to $\bar{\e}\in \bar{E}$ will be denoted by
$\bar{e}_\e$.

We define a homomorphism
$\hat{\phi}: \La \to \La'$
by specifying the images of generators.
Specifically,
\begin{alignat*}{2}
\hat{\phi}(e_i)&=\bar{e}_i-\bar{e}_n+e_p && \quad \text{for} \ i \in [n-1],  \\
\hat{\phi}(e_n)&=e_p,  && \quad \text{and}\hfill  \\
\hat{\phi}(e_\e)&=\bar{e}_\e && \quad \text{for} \ \e \in E_0. 
\end{alignat*}

\begin{lemma}The map $\hat{\phi}:\La \to \La'$ is an isomorphism.
\end{lemma}

\begin{proof} Keeping in mind that $\bar{e}_1=\bar{e}_{\e_0}$ in $\La'$, we see that $\hat{\phi}$ has a well-defined inverse in degree one
given by 
\begin{alignat*}{2}
\bar{e}_i & \mapsto e_i-e_1+e_{\e_0} && \quad \text{for $1\leq i\leq n$,} \\  
\bar{e}_\e & \mapsto e_\e && \quad \text{for $\e \in E_0$, and} \\
e_p& \mapsto e_n. &&
\end{alignat*} 
It follows that $\hat{\phi}$ is an isomorphism.
\end{proof}

\begin{lemma} $\hat{\phi}(I) \subseteq I'$.
\end{lemma}

\begin{proof}
If $\{\e_1,\ldots,\e_q\}$ is a circuit of $M_0$, then $\hat{\phi}(\partial
e_{\e_1}\cdots e_{\e_q})=\partial \bar{e}_{\e_1}\cdots \bar{e}_{\e_q}.$
With the observation that $\hat{\phi}(e_i-e_{i+1})=\bar{e}_i-\bar{e}_{i-1}$
for $1\leq i \leq n-2$, and also for $i=n-1$, we see that 
\begin{equation*}
\begin{split}
\hat{\phi}(\partial e_1\cdots e_n)&=\hat{\phi}((e_1-e_2)(e_2-e_3)\cdots
(e_{n-1}-e_n)) \\
& = (\bar{e}_1-\bar{e}_2)(\bar{e}_2-\bar{e}_3)\cdots(\bar{e}_{n-1}-\bar{e}_n)
\\& = \partial\bar{e}_1\cdots \bar{e}_n.\\
\end{split}
\end{equation*}
Referring to the definitions of $M$ and $M'$, we see that these computations
suffice to prove the lemma.
\end{proof}

\begin{corollary} $\hat{\phi}: \La \to \La'$ induces a surjection $\phi:
A \to A'$.
\label{onto}
\end{corollary}

\end{section}

\begin{section}{The Tutte polynomials}

By the end of this section we will have proved Theorem \ref{main}. The final
ingredient is the computation of Tutte polynomials.
The Tutte polynomial $T_M(x,y)$ is defined recursively as follows.
$M\setminus e$ and $M/e$ refer to the deletion and contraction of $M$
relative to $e$. 
\begin{enumerate}
\item $T_M(x,y)=x$ if $M$ is an isthmus; $T_M(x,y)=y$ if $M$ is a loop.
\item $T_M(x,y)=T_e(x,y)T_{M\setminus e}(x,y)$ if $e$ is a loop or isthmus
in $M$.
\item $T_M(x,y)=T_{M\setminus e}(x,y)+T_{M/e}(x,y)$ otherwise.
\end{enumerate}
These properties uniquely determine a polynomial
$T_M(x,y)$ which is a matroid-isomorphism invariant of $M$.

We will use the following standard property of Tutte polynomials.
\begin{lemma} $T_{M\oplus M'}(x,y)=T_M(x,y)T_{M'}(x,y).$
\label{product}
\end{lemma}

The characteristic polynomial $\chi_M(t)$ of $M$ may be defined by
$$\chi_M(t)=T_M(1-t,0).$$ The following result of \cite{OS1} was the initial
cause for interest in $A(M)$ among combinatorialists. We will use it to
show that $\phi$ is injective.

\begin{theorem} The Hilbert series $$\sum_{p=0}^\infty \dim(A^p)t^p$$ of $A=A(M) 
$
is equal to $t^r\chi_M(-t^{-1})$, where $r=\rk(M)$.
\label{hilbert}
\end{theorem}

In fact the $OS$ algebra is isomorphic to the Whitney homology of the
lattice of flats of $L$, equipped with a natural product \cite{Bj}.

The next lemma is easy to prove by induction on $n$.

\begin{lemma} For any $n\geq 2$, $T_{C_n}(x,y)=\sum_{i=1}^{n-1}x^i +y.$
\label{cycle}
\end{lemma}
 
Lemma \ref{cycle} and Theorem \ref{tutte} may be deduced from more general results proved in Section 6 of \cite{Bry71}. We include the proof of \ref{tutte} here for the reader's convenience. Let $M$ and $M'$  be the matroids of the preceding section. 

\begin{theorem} Let $n\geq 2$, Then 
$$T_M(x,y)= \bigl(\sum_{i=1}^{n-1}x^i +y\bigr)T_{M_0}(x,y),  \ \text{and}$$ 
$$T_{M'}(x,y)=\bigl(\sum_{i=1}^{n-1}x^i\bigr)T_{M_0}(x,y) + xyT_{M_0/\e_0}(x,y).$$
\label{tutte}
\end{theorem}

\begin{proof}
The first formula is a consequence of Lemmas \ref{cycle} and
\ref{product}. To prove the second assertion,  we establish a recursive
formula for the Tutte polynomial of $P^n_{\e_0}$. Assume $n\geq 3$, and
apply property (iii) above to a point of $C_n$ other than 1. The deletion
is the direct sum of $M_0$ with $n-2$ isthmuses, and the contraction is
$P^{n-1}_{\e_0}$. Thus we have $$T_{P^n_{\e_0}}(x,y)=x^{n-2}T_{M_0}(x,y)+
T_{P^{n-1}_{\e_0}}(x,y).$$

Now consider the case $n=2$. Deleting the point $\bar{2}$ yields $M_0$,
while contracting $\bar{2}$ yields the direct sum of $M_0/\e_0$ with a
loop. Thus $$T_{P^2_{\e_0}}(x,y)=T_{M_0}(x,y)+yT_{M_0/\e_0}(x,y).$$
Then one can prove inductively that
$$T_{P^n_{\e_0}}(x,y)=\bigl(\sum_{i=0}^{n-2}x^i\bigr)T_{M_0}(x,y)
+yT_{M_0/\e_0}.$$
Since $M'$ is the direct sum of $P^n_{\e_0}$ with an isthmus, right-hand
side of this formula is multiplied by $x$ to obtain $T_{M'}(x,y)$.
\end{proof}

\begin{corollary} $\chi_M(t)=\chi_{M'}(t)$
\label{whitney}
\end{corollary}

\begin{proof} The two formulas in Theorem \ref{tutte} yield the same
expression upon setting $y=0$. The assertion then follows from the
definition of $\chi_M(t)$ above.
\end{proof}

\begin{corollary} The map $\phi: A \to A'$ is an isomorphism.
\label{isom}
\end{corollary}

\begin{proof} According to Theorem \ref{hilbert}, the last corollary implies
$\dim A^p =\dim (A')^p $. Since $\phi$ is surjective by \ref{onto}, and all
spaces are finite-dimensional, $\phi$ must be an isomorphism.
\end{proof}

In case $n=3$ and $M_0=C_3$, 
the map $\phi$ is a modified version of the isomorphism discovered by 
L.~Rose and H.~Terao \cite[Example 3.77]{OT} for the rank three 
truncations of $M_3$ and $M'_3$. 

With the next result, we complete the proof of Theorem \ref{main}.
\begin{corollary} If $M_0$ is connected, 
then $T_M(x,y)\not = T_{M'}(x,y)$.
\label{isthmus}
\end{corollary}

\begin{proof} Assume $T_M(x,y)=T_{M'}(x,y)$. By Theorem \ref{tutte} 
this implies $$T_{M_0}(x,y)=xT_{M_0/\e_0}(x,y).$$ 
By hypothesis $\e_0$ is not an isthmus. Deleting and contracting along
$\e_0$, and evaluating at $(x,y)=(1,1)$, we obtain $$T_{M_0\setminus
\e_0}(1,1)+T_{M_0/\e_0}(1,1)=T_{M_0/\e_0}(1,1),$$ which implies
$T_{M_0\setminus \e_0}(1,1)=0$. Coefficients of Tutte polynomials are
non-negative, so this implies $T_{M_0\setminus \e_0}(x,y)=0$, which is not
possible.
\end{proof}

The proof of the last corollary uses only the fact that $\e_0$ is not an
isthmus. Thus Theorem \ref{main} remains true for any simple matroid $M_0$ which
is not the uniform matroid of rank $m$ on $m$ points (realized by the boolean
arrangement of coordinate hyperplanes), in which every point is an isthmus.

\begin{remark}
The proof of Corollary \ref{isthmus} specializes, upon setting
$(x,y)=(1-t,0)$,
to a proof of the result of H.~Crapo that a connected matroid has nonzero beta
invariant \cite{Wh1}.
\end{remark}
\end{section}

\begin{section}{Concluding remarks}

We start this section with a proof of Corollary \ref{cor}. Let $G_m$ be the
graph with vertex set $\Z_{2m}$ and edges $\{i,i+1\}$ for
$1\leq i<2m-1$ and
$\{0,i\}$ for $1\leq i< 2m$. Then $G_m$ has $2m$ vertices and $4m-3$ edges. The graph $G_4$ is illustrated in Figure \ref{wheel}.

%\psdraft
\begin{figure}
\begin{center}
\epsfig{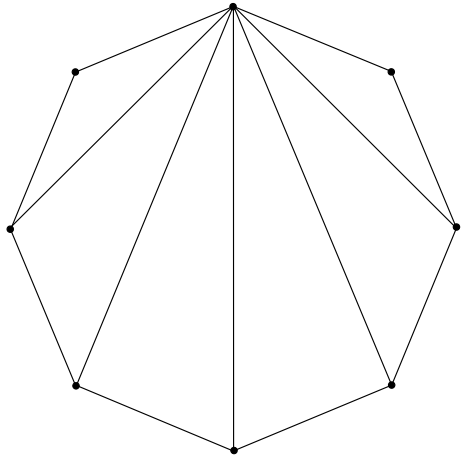}
\caption{The graph $G_4$.} \label{wheel}
\end{center}
\end{figure}
%\psfull

\begin{theorem}Let $n>2m+1$. Then
the parallel connections of $G_m$ with $C_n$ along the edges $\{0,i\}$
of $G_m$
result in mutually non-isomorphic graphs, for $m\leq i \leq 2m-1$.
\label{parallel}
\end{theorem}

\begin{proof} Fix $i$ in the specified range.
Then the parallel connection $P(C_n,G_m)$ along $\{1,i\}$ has longest
circuit of
length $(n-1)+i$. The assertion follows.
\end{proof}

\begin{proofof}{Corollary \ref{cor}}
Let $M_0$ be the polygon matroid of $G_m$. The proof of Theorem
\ref{parallel} actually
shows that the parallel connections $P^n_i$ of $M_0$ with $C_n$ along $\{1,i\}$ yield $m$
non-isomorphic matroids  as $i$ ranges from $m$ to $2m$. The same holds
true when they are extended by an isthmus, resulting in $m$ non-isomorphic
matroids $M'_{n,i}$. But the proof of Corollary \ref{isom} did not depend
on the choice of $\e_0$. So the $OS$ algebra of $M'_{n,i}$ is isomorphic to
the $OS$ algebra of $M_n=C_n\oplus M_0$ independent of $i$. 
This completes the proof of \ref{cor}.
\end{proofof}

We close with some topological considerations. 
We will see that part of Theorem \ref{main}, in the case that $M_0$ is
realizable over \C, is a consequence of a general topological equivalence. This equivalence follows from a well-known relationship
between the complements of a central arrangement and its projective image.
The proof is quite trivial, but requires us to introduce
explicit realizations, with apologies for the cumbersome
notation.
We will need a few easy facts about hyperplane complements, which may be
found in \cite{OT}.

Let $\A=\{H_1,\ldots,H_n\}$ be an arrangement of affine hyperplanes in
$\C^r$. Let $\phi_i: \C^r \to \C$ be a linear polynomial function with 
$H_i=\{x\in\C^r \ | \ \phi_i(x)=0\}$. The {\em defining polynomial} of
\A\ is the product
$Q(\A)=\Pi_{i=1}^n \phi_i$. If all of the $\phi_i$ are homogeneous linear
forms, \A\ is said to be a {\em central} arrangement. 

Recall $C(\A)$
denotes the complement of $\bigcup \A$ in $\C^r$. The connection
between central arrangements in $\C^r$ and affine arrangements in
$\C^{r-1}$ goes as follows. Assume \A\ is central. 
Change variables so that $\phi_1(x)=x_1$,
and write $Q(\A)=x_1\hat{Q}(x_1,\ldots,x_r)$.
Consider $(x_2,\ldots,x_r)$ to be coordinates on $\C^{r-1}$. 
Then let $d\A$ denote the affine arrangement in $\C^{r-1}$
with defining polynomial
$\hat{Q}(1,x_2,\ldots,x_r)$.

\begin{lemma} $C(\A)$ is diffeomorphic to $\C^* \times C(d\A)$.
\label{affine}
\end{lemma}

Suppose $\A_0$ and $\A_1$ are affine arrangements with defining polynomials
$Q_0(x)$ and $Q_1(y)$ in disjoint sets of variables 
$x=(x_1,\ldots,x_{r_0})$ and $y=(y_1,\ldots,y_{r_1})$. 
Let $\A_0\oplus \A_1$ be the arrangement in
$\C^{r_0+r_1}$ with defining polynomial $Q_0(x)Q_1(y)$.

\begin{lemma} $C(\A_0\oplus \A_1)$ is diffeomorphic to $C(\A_0)\times
C(\A_1)$.
\end{lemma}

If $\A_0$ and $\A_1$ are central arrangements, with underlying matroids
$M_0$ and $M_1$, then $\A_0\oplus \A_1$ is a realization of the direct sum
$M_0\oplus M_1$.

Now let $\A_0$ and $\A_1$ be arbitrary central arrangements, with
underlying matroids $M_0$ and $M_1$. 
To realize the parallel connection $P(M_0,M_1)$ change coordinates so that
the hyperplane $x_1=0$ appears in both $\A_0$ and $\A_1$. These will be the
hyperplanes that get identified in the parallel connection. Write
$$Q_1(\A_1)=x_1\hat{Q}_1(x_1,\ldots,x_{r_1}).$$ 
With $(x_1,\ldots,x_{r_0},y_2,\ldots, y_{r_1})$ as
coordinates in $\C^{r_0+r_1-1}$, the parallel connection $P(\A_0,\A_1)$ is 
the arrangement in $\C^{r_0+r_1-1}$ defined by the polynomial
$$Q_0(x_1,\ldots,x_{r_0})\hat{Q}_1(x_1,y_2,\ldots,y_{r_1})$$

We are now prepared to prove Theorem \ref{homo}.
Let ${\cal S}$ denote the arrangement in $\C^1$ with defining polynomial
$x$. So ${\cal S}$ has as underlying matroid the isthmus $S$, and $C({\cal
S})=\C^*$.

\begin{proofof}{Theorem \ref{homo}}
Write $Q(\A_0)=x_1\hat{Q}_0(x_1,\ldots,x_{r_0})$.
Following the recipe given above for dehomogenizing an arrangement, and
using the given defining polynomial for $P(\A_0,\A_1)$, 
we see that the affine arrangement
$\bar{P(\A_0,\A_1)}$ has defining polynomial
$$\hat{Q}_0(1,x_2,\ldots,x_{r_0})\hat{Q}_1(1,y_2,\ldots,y_{r_1}),$$
which is precisely the defining polynomial of $d\A_0\oplus
d\A_1.$
By the preceding lemmas we have 
\begin{equation*}
\begin{split}
C({\cal S}\oplus P(\A_0,\A_1)) & \cong C(\S) \times C(P(\A_0,\A_1) \\
& \cong \C^* \times \C^* \times C(\bar{P(\A_0,\A_1)}) \\
& \cong \C^* \times \C^* \times C(d\A_0) \times C(d\A_1).\\
\end{split}
\end{equation*}

On the other hand,
$$C(\A_0\oplus \A_1)\cong C(\A_0)\times
C(\A_1)\cong \C^* \times C(d\A_0)\times \C^* \times C(d\A_1).$$
This proves the result.
\end{proofof}

Returning to Theorem \ref{main}, if we assume the matroid 
$M_0$ is realizable over \C,
we can take such a realization for $\A_0$, and any general position
arrangement of $n$ hyperplanes in $\C^{n-1}$ for $\A_1$. Then 
Theorem \ref{homo} and Theorem \ref{OS} together imply that the 
$OS$ algebras $A(M_n)$ and $A(M'_n)$ over \C\ are isomorphic.

The arrangements constructed
in \cite{F3} are generic 3-dimensional sections of the arrangements
of Theorem \ref{homo}, with the seed arrangements $\A_0$ and
$\A_1$ both of rank two. The fact that their fundamental groups are isomorphic is then an immediate consequence of the Lefschetz Hyperplane Theorem. This theorem does not imply that the sections are homotopy equivalent; this is proved in \cite{F3} by constructing an explicit isomorphism of canonical presentations of the fundamental groups, using Tietze transformations only of type I and II. We do not know if the diffeomorphic arrangements constructed in Theorem \ref{homo} will in general have homotopy equivalent generic 3-dimensional sections.

\medskip
The question whether arrangements with different combinatorial structure 
could have homotopy equivalent complements was originally restricted to
central arrangements because Lemma \ref{affine} provides trivial 
counter-examples in the affine case.
The constructions presented in this paper are now seen from the  proof of
Theorem \ref{homo} to arise again from Lemma \ref{affine}. So the examples of
\cite{F3} also come about in some sense from Theorem \ref{homo}.  
We feel compelled to again narrow the problem to rule out these other, 
not quite so
trivial counter-examples.

\begin{conjecture} For central arrangements whose underlying matroid is
connected, the homotopy type of the complement determines the underlying
matroid.
\end{conjecture}

\end{section}

\begin{ack} This paper was completed during the second author's visit to
the Mathematical Sciences Research Institute in Berkeley 
in the early spring of 1997. He is grateful to the staff and
organizers for their support. He also thanks Richard Stanley for 
the last step in the proof of Corollary \ref{isthmus}. 
We are both thankful to Terry
Blows for his work in administering the REU program at Northern Arizona
University during which this research was done. Finally we thank the referee for pointing us toward the reference \cite{Bry71}.
\end{ack}

%\bibliography{biblio}
%\bibliographystyle{plain}

\end{document}